\documentclass[12pt]{article}
\usepackage{amsmath}
\usepackage{graphicx,psfrag,epsf}
\usepackage{enumerate}
\usepackage{natbib}
\usepackage{url} 

\usepackage{graphicx}
\usepackage{slashbox}
\usepackage{amssymb}
\usepackage[english]{babel}
\usepackage{amssymb,amsmath,multirow,array}
\usepackage[english]{babel}
\newtheorem{theorem}{Theorem}
\newtheorem{lemma}{Lemma}

\newcommand{\blind}{0}

\date{}
\begin{document}


\def\spacingset#1{\renewcommand{\baselinestretch}%
{#1}\small\normalsize} \spacingset{1}


  \title{Asymptotic properties in the Probit-Zero-inflated Binomial regression model}
 \author{DIOP Aba\footnote{Corresponding author}\\
 \scriptsize{Equipe de Recherche en Statistique et Modèles Aléatoires (ERESMA),}\\
   \scriptsize{Université Alioune Diop, Bambey, Sénégal}\\
   \scriptsize{(aba.diop@uadb.edu.sn)}\\
 BA Demba Bocar   \\
 \scriptsize{Département de Maths, Université Iba Der Thiam, Thiès, Sénégal} \\
 \scriptsize{(dbba@univ-thies.sn)} \\
 LO Fatimata\\
 \scriptsize{Département de Maths, Université Iba Der Thiam, Thiès, Sénégal} \\
 \scriptsize{(fatymalo@hotmail.com)} }
  
\maketitle
\if0\blind

\bigskip
\begin{abstract}
Zero-inflated regression models have had wide application recently and have provenuseful in modeling data with many zeros. Zero-inflated Binomial (ZIB) regression model is an extension of the ordinary binomial distribution that takes into account the excess of zeros. In comparing the probit model to the logistic model, many authors believe that there is little theoretical justification in choosing one formulation over the other in most circumstances involving binary responses. The logit model is considered to be computationally simpler but it is based on a more restrictive assumption of error independence, although many other generalizations have dealt with that assumption as well. By contrast, the probit model   assumes that random errors have a multivariate normal distribution. This assumption makes the probit model attractive  because the normal distribution provides a good approximation to many other distributions. In this paper, we develop a maximum likelihood estimation procedure for the parameters of a zero-inflated Binomial regression model with probit link function for both component of the model. We establish the existency, consistency and asymptotic normality of the proposed estimator.
\end{abstract}

\noindent%

\spacingset{1.45} 
\section{Introduction}
\label{sec:intro}
Zero-inflated regression models mix a degenerate distribution with point mass of one at zero with a simple regression model based on a standard distribution. Zero-infation can happen both with discrete data and continuous data (see \cite{lambert92}, \cite{h2000} and \cite{famoy06}). Here, we  concentrate on the discrete case, specialy on Zero-inflated Binomial (ZIB) regression model. \cite{h2000} considered Zero-Inflated Poisson (ZIP) models and also extended Lambert's methodology to an upper bounded count situation, thereby obtaining a ZIB model. He also incorporated random intercepts into these models and used them to model insect counts in a horticulture experiment. The ZIB regression model assumes the data are from two states: a true zero state where the response always equals zero and a non zero state where the response follows a Binomial distribution. Suppose we have the response vector $y=(y_1,...,y_n)'$ where $y_i$ is the observed value of a random variable $Y_i$ and $n\in\mathbb{N}^{*}$ is the sample size. In the ZIB regression model we assume
\begin{equation}\label{sec1eq1}
Y_i \sim \begin{cases}
0 & \text{with probability}~~ p_i \\
\mathcal{B}(n_i,\pi_i)  & \text{with probability}~~ 1-p_i
\end{cases}
\end{equation}
where $y_i$ is the number of successes out of $n_i$ trials and $\pi_i$ is the probability of success for the subject $i$. Again, the parameters $p=(p_1,\ldots, p_n)'$ and $\pi_i=(\pi_1,\ldots, \pi_n)'$ are modeled through \textit{logit} link generalized linear models as:
\begin{equation}\label{sec1eq2}
\begin{split}
\log\left[ \frac{p_i}{1-{p_i}}\right]  &=\theta_1+\theta_2\mathbf z_{i2}+\ldots +\theta_p \mathbf z_{iq}:= \theta'\mathbf z_i  \\
\log\left[ \frac{\pi_i}{1-{\pi_i}}\right] &= \beta_1+\beta_2\mathbf x_{i2}+\ldots +\beta_p \mathbf x_{ip}:= \beta' \mathbf x_i
\end{split}
\end{equation}
Here, $\beta$ and $\theta$ are regression parameters and $\mathbf X$ and $\mathbf Z$ are corresponding design matrices. \cite{alpha17} studied the asymptotics of the model \ref{sec1eq1} by considering the logit link function for the two components of the model (binomial and zero inflation). They established the consistency and asymptotic normality of the maximum likelihood estimator. \cite{v2000} applied the ZIB model estimated by a quasi-likelihood method to a biological control assay. This model was initially built to model overdispersion generated by individual variability in the probability of success.\\
For other types of discrete outcomes, such as binary, multinomial or ordinal, various single value inflated models were developed, including: zero-inflated Bernoulli model (\cite{ddd16}), zero-inflated binomial model (\cite{h2000} and \cite{v2000}), zero-inflated ordered probit model (\cite{ha2007}), baseline or zero inflated multinomial logit model (\cite{bag16}, \cite{alpha17b}), and middle category inflated ordered model (\cite{bag12}). Similar extension has been made to incorporate inflation other than zero for multinomial or ordinal outcomes (\cite{sweeney14}). The use of noncanonical link functions is not prohibited by the fact that they are more computationally complex. \cite{czado2000} argued that, in some applications, the overall fit of the model as measured by the p-value of the  goodness-of-fit statistics can be improved significantly by the use of a noncanonical link. \\
We use the probit link function for both  binomial and excess of zero components for the ZIB regression model \ref{sec1eq1}. The aim is to evaluate the misspecification problem by mistakenly choosing the appropriate link function. We do not introduce the probit  model as a rival to the logistic model, but rather as an alternative, in particular in the case of excess of zero in the sample. Experience shows that in most situations the two approaches produce similar results although some differences exist. This similarity is not necessarily sustained when multivariate responses are used. Further research is needed to investigate the advantages or disadvantages in using one model over the other in data mining applications. \\

The use of the probit regression model dates back to \cite{bliss34}. Many other authors have used the probit model in other applications with success; for example, \cite{sha09} compared the two models for estimating the strength of gear teeth. \cite{raz13} used probit link function for data mining applications. \cite{fin71} applied probit analysis in toxicological experiments. The probit model has also found popularity in economics. \cite{cramer03} provides a survey of the early origins of the model. \\

The rest of this paper is organized as follows. In Section \ref{sec:model}, we describe the problem of ZIB regression model with probit link function, we propose an estimation method adapted to this setting. In Section \ref{sec:asymptotic} we establish the asymptotic propreties of the proposed estimator under some regularity conditions. A discussion and some perspectives are given in Section \ref{sec:discussion}.

\section{ZIB regression model with Probit link function}
\label{sec:model}
\subsection{The model set-up and estimation procedure}
Let $(Y_1, \mathbf X_1, \mathbf W_1),\ldots, (Y_n, \mathbf X_n, \mathbf W_n)$ be independent and identically distributed copies of the random vector $(Y, \mathbf X, \mathbf W)$ defined on the probability space $(\Omega, \mathcal A, \mathbb P)$. For every individual $i=1, \ldots,n$, $Y_i$ is the number of successes out of $n_i$ trials and $\pi_i$ is the probability of success for the subject $i$. Let $\mathbf X_i=(1, X_{i2}, \ldots, X_{ip})^{\top}$ and $\mathbf W_i=(1, Z_{i2}, \ldots, W_{iq})'$ be random vectors of predictors or covariates (both categorical and continuous predictors are allowed). We shall assume in the following that the $\mathbf X_i$'s are related to the zero inflation, while the $\mathbf W_i$'s are related to the binomial component. $\mathbf X_i$ and $\mathbf W_i$ are allowed to share some components.\\
The ZIB regression model with \textit{probit} link function is given by
\begin{equation}\label{sec2eq1}
Y_i \sim \begin{cases}
0 & \text{with probability}~~ p_i \\
\mathcal{B}(n_i,\pi_i)  & \text{with probability}~~ 1-p_i
\end{cases}
\end{equation}
where the parameters $p=(p_1,\ldots, p_n)'$ and $\pi_i=(\pi_1,\ldots, \pi_n)'$ are modeled through \textit{probit} link generalized linear models as:
\begin{equation}
\begin{split}
\text{\textit{probit}}(p_i)  &= F(\beta_1+\beta_2\mathbf x_{i2}+\ldots +\beta_p \mathbf x_{iq}):= F(\beta'\mathbf x_i)  \\
\text{\textit{probit}}(\pi_i) &= F(\mu_1+\mu_2\mathbf w_{i2}+\ldots +\mu_p \mathbf w_{ip}):= F(\mu' \mathbf w_i)
\end{split}
\end{equation}
where $F$ is the cumulative distribution function of the standard normal distribution.\\\\
The model (\ref{sec2eq1}) can be rewritten as follows
\begin{equation}\label{sec2eq2}
\mathbb{P}(Y=y_i)=
 \begin{cases} 
\pi_i +(1- \pi_i)(1-p_i)^{n_i} &  \text{if}~~  y_i= 0\\
(1 - \pi_i)C_{n_i}^{y_i}p_i^{y_i}(1-p_i)^{n_i-y_i} &  \text{if}~~ y_i \in \{1, \ldots, n\}
\end{cases}
\end{equation}
~\\
Let $\theta:=(\beta',\mu')'$ denote the unknown $k$-dimensional ($k=p+q$) parameter in the conditional distribution of $Y$ given $\mathbf X_i$ and $\mathbf W_i$. Now, the likelihood for $\theta$ from the independent sample $(Y_i, \mathbf X_i, \mathbf W_i)$ $(i=1,\ldots,n)$ is as follows:
\begin{equation}\label{sec2eq3}
\begin{split}
L_n(\theta) &= \prod_{i=1}^n \mathbb{P}(Y_i=k_i | \mathbf X_i = \mathbf x_i,\mathbf W_i = \mathbf w_i) \\
 &= \prod_{i=1}^{n} (F(\mu'\mathbf w_i) + (1-F(\mu'\mathbf w_i))(1- F(\beta'\mathbf x_i))^{n_i})^{k_i} \\
 &\times ((1- F(\mu'\mathbf w_i))C_{n_i}^{Y_i}F(\beta'\mathbf x_i)(1-F(\beta'\mathbf x_i))^{n_i-Y_i})^{1-k_i}
\end{split}
\end{equation}
where $k_i = \mathbf 1_{ \lbrace Y_i=0\rbrace}$.\\
We define the maximum likelihood estimator $\widehat \theta_n:=(\widehat \beta'_{n},\widehat \mu'_{n})'$ of $\theta$ as the solution (if it exists) of the $k$-dimensional score equation
%
\begin{equation}\label{sec2eq3}
\dot{l}_n(\theta)=\frac{\partial l_n(\theta)}{\partial\theta}=0,
\end{equation}
where $l_n(\theta):=\log L_n(\theta)$ is the log-likelihood function. In the following, we shall be interested in the asymptotic properties of the maximum likelihood estimator $\widehat \theta_n$. We will however obtain consistency
and asymptotic normality results for the whole $\widehat \theta_n$. Before proceeding, we need to set some further notations.
\subsection{Some further notations}
Define first the $(p\times n)$ and $(q\times n)$ matrices
\begin{equation*}
\mathbb X=
\begin{pmatrix}
1& 1 & \dots & 1\\
X_{12} & X_{22} &\dots & X_{n2}\\
\vdots & \vdots & \ddots & \vdots \\
X_{1p} & X_{2p} &\ldots & X_{np}
 \end{pmatrix}
 \quad
  \text{and} 
  \quad  
 \mathbb W =
 \begin{pmatrix}
1& 1 & \dots & 1\\
W_{12} & W_{22} &\dots & W_{n2}\\
\vdots & \vdots & \ddots & \vdots \\
W_{1p} & W_{2p} &\ldots & W_{np}
 \end{pmatrix}
\end{equation*}
and let $\mathbb Z$ be the $(k\times 2n)$ block-matrix defined as
\begin{equation*}
\mathbb Z =
\begin{bmatrix}
\mathbb X & 0_{pn}\\
0_{qn} & \mathbb W \\
\end{bmatrix}
\end{equation*}
where $0_{ab}$ denotes the $(a\times b)$ matrix whose components are all equal to zero (for any positive integer values $a$, $b$). Let $k_i(\theta) =F(\mu'\mathbf w_i)(h_i(\beta))^{-n_i+1}+(1-F(\mu'\mathbf w_i)h_i(\beta)  $, $\forall i =1,\ldots, n$ where $h_i(\beta) = 1 - F(\beta'\mathbf x_i)$. Let also $C(\theta) = (C_j(\theta))_{1\leq j\leq 2n}$ be the 2n-dimensional column vector defined as
\begin{equation*}
C(\theta) = (A_1(\theta),\ldots, A_n(\theta);B_1(\theta),\ldots,B_n(\theta))', \quad i= 1, \ldots, n
\end{equation*}
where
\begin{equation*}
\begin{split}
A_i(\theta) &=  \frac{k_i\mathbf w_i f(\mu'\mathbf w_i)h_i(\beta)^{-n_i+1}F(\beta'\mathbf x_i)^{n_i}}{k_i(\theta)} -(1-k_i)\frac{\mathbf w_i f(\mu'\mathbf w_i)}{1-F(\mu'\mathbf w_i)} \\
B_i(\theta) &= \frac{-k_i n_i \mathbf x_i f(\beta'\mathbf x_i)}{F(\mu'\mathbf w_i)k_i(\theta)} + (1-k_i)\frac{\mathbf x_i f(\beta'\mathbf x_i)[\mathbb z_i-n_i F(\beta'\mathbf x_i)]}{F(\beta'\mathbf x_i)h_i(\beta)}     
\end{split}
\end{equation*}
where $f(.)$ is the density function associate to $F(.)$.\\ \\
Then, simple algebra shows that the score equation can be rewritten as
\begin{equation*}
\dot{l}_n(\theta) = ZC(\theta) = 0.
\end{equation*}
If $H= (H_{ij})_{1\leq i \leq a,1\leq j\leq b}$ denotes some $(a\times b)$ matrix, we will denote by $H_{\bullet j}$ it's $j$-th column $(j = 1,\ldots, b)$ that is, $H_{\bullet j} = (H_{1j},\ldots, H_{aj})'$. Then, it will be useful to rewrite the score vector as
\begin{equation*}
\dot{l}_n(\theta) = \sum_{i=1}^{2n} \mathbb Z{\bullet j}C_j(\theta).
\end{equation*}
We shall further note $\ddot{l}_n(\theta)$ the $(k\times k)$ matrix of second derivatives of $l_n(\theta)$, that is $\ddot{l}_n(\theta) = \frac{\partial^2 l_n(\theta)}{\partial\theta \partial\theta'}$. Let $\mathbf D(\theta)=(\mathbf D_{ij}(\theta))_{1\leq i,j \leq 2n}$ be the $(2n\times 2n)$ block matrix defined as
\begin{equation*}
\mathbf{D}(\theta) = 
\begin{bmatrix}
\mathbf D_1(\theta) & \mathbf D_3(\theta)\\
\mathbf D_3(\theta) & \mathbf D_2(\theta)\\
\end{bmatrix}
\end{equation*}
where $\mathbf D_1(\theta), \mathbf D_2(\theta)$ and $\mathbf D_3(\theta)$ are $(n\times n)$ diagonal matrices, with $i$-th diagonal elements $(i= 1,\ldots, n)$ respectively given by
\begin{equation*}
\begin{split}
\mathbf D_{1,ii}(\theta) &= \frac{k_{i}(\theta) w_{i}^{2} h_i(\beta)^{-n_i+1}F(\beta'\mathbf x_i)[f'(\mu'\mathbf w_i)k_i(\theta)-f^{(2)}(\beta\mathbf x_i)(h_i(\beta)^{-n_i+1}-h_i(\beta)]}{k_{i}^{2}(\theta)} \\
&- \mathbf w_{i}^{2}(1-k_{i}(\theta))\frac{[f'(\mu\mathbf w_i)(1-F(\mu'\mathbf w_i))+f^{(2)}(\mu'\mathbf w_i)]}{(1-F(\mu'\mathbf w_i))^2}  \\
\mathbf D_{2,ii}(\theta) &= \frac{-k_{i}(\theta) n_i \mathbf x_i}{F(\mu'\mathbf w_i}[x_i f(\beta\mathbf x_i)k_i(\theta)-(n_i+1)x^{2}_{i} f^{(2)}(\beta'\mathbf x_i)(F(\mu'\mathbf w_i)h_i(\beta)^{n_i}-(1-F(\mu'\mathbf w_i))] \\
&+ (1-k_{i}(\theta))\frac{[\mathbf x_{i}^{2} F(\beta'\mathbf x_i)h_i(\beta)(\mathbf z_i f'(\beta'\mathbf x_i)-n_i f'(\beta'\mathbf x_i)F(\beta'\mathbf x_i)}{(F(\beta'\mathbf x_i)h_i(\beta))^2} \\
&- (1-k_{i}(\theta))\frac{n_i f^{(2)}(\beta'\mathbf x_i)) + \mathbf x_i f(\beta'\mathbf x_i)(\mathbf z_i-n_i F(\beta'\mathbf x_i)(x_i f(\beta'\mathbf x_i)h_i(\beta) -x_i f(\beta'\mathbf x_i)F(\beta'\mathbf x_i))}{(F(\beta'\mathbf x_i)h_i(\beta))^2}\\
\mathbf D_{3,ii}(\theta) &= \frac{-k_i \mathbf w_i f(\mu'\mathbf w_i)k_i(\theta)[(-n_i+1)\mathbf x_i h_i(\beta)^{-n_i}F(\beta'\mathbf x_i)+\mathbf x_i f(\beta'\mathbf x_i)h_i(\beta)^{-n_i+1})}{k_i^2(\theta)}\\
&- \frac{k'_i(\theta)h_i(\beta)^{-n_i+1} F(\beta'\mathbf x_i)]}{k_i^2(\theta)}
\end{split}
\end{equation*}
Then, some algebra shows that $\ddot{l}_n$ can be expressed as
\begin{equation*}
\ddot{l}_n(\theta) = -\mathbb{Z}\mathbb{D}(\theta)\mathbb{Z}'.
\end{equation*}
Note that the size of $C(\theta)$,$\mathbb Z$ and $\mathbf D(\theta)$ depends on $n$. However, in order to simplify
notations, $n$ will not be used as a lower index for these vector and matrices. In the next section, we turn to the asymptotic theory for the proposed estimator.
\section{Asymptotic theory}
\label{sec:asymptotic}
We first state some regularity conditions that will be needed to study the asymptotic theory. 
\begin{description}
\item[C1] The covariates are bounded that is, there exist compact sets $G\subset\mathbb R^p$ and $H\subset\mathbb R^q$ such that $\mathbf X_i\in G$ and $\mathbf W_i\in H$ for every $i=1,2,\ldots$ For every $i=1,2,\ldots$, $j=2,\ldots,p$, $k=2,\ldots,q$, $\mbox{var}[X_{ij}]>0$ and $\mbox{var}[W_{ik}]>0$. For every $i=1,2,\ldots$, the $X_{ij}$
$(j=1,\ldots,p)$ are linearly independent, and the $W_{ik}$ $(k=1,\ldots,q)$ are linearly independent.
\item[C2] Let $(\beta'_{0}, \mu'_{0})'$ denote the true parameter value. $\beta_0$ and $\mu_0$ lie in the interior of known compact sets $\mathcal G \subset\mathbb R^p$ and $\mathcal H\subset\mathbb R^q$ respectively.
\item[C3] The Hessian matrix $\ddot{l}_n(\theta)$ is negative definite and of full rank, for every $n=1,2,\ldots$ Let $\lambda_n$ and $\Lambda _n$ be respectively the smallest and largest eigenvalues of $\mathbb Z \mathbf D(\theta_0) \mathbb Z'$. There exists a finite positive constant $c_2$ such that $\Lambda_n\slash\lambda_n<c_2$ for every $n=1,2,\ldots$
\item[C4] Let $A$ denote a $k\times k$ matrix such that $\|A\| = \max_{1\leq i\leq k}|\lambda_i|$.
\item[C5] The function $f$ is continuous and twice differentiable.
\end{description}
Now we establish rigorously the existence, consistency and asymptotic normality of the maximum likelihood estimator $\widehat{\theta}_n$ in model (\ref{sec2eq1}). We prove the following results:
%
%
\begin{theorem}[Existence and consistency]\label{th1}
Under the conditions \textbf{C1}-\textbf{C5}, The maximum likelihood estimator $\widehat \theta_n$ exists almost surely as $n\rightarrow \infty $ and converges almost surely to $\theta_0$, if and only if $\lambda_n$ tends to infinity as $n\rightarrow \infty$.
\end{theorem}
\noindent
\textbf{Proof of Theorem \ref{th1}} The following lemma essentially provides an intermediate technical result (see \cite{alpha17} for its proof).
\begin{lemma}\label{lem1}
Let $\phi_n(\theta): \mathbb{R}^{k}\longrightarrow \mathbb{R}^{k}$ be defined as
\begin{equation*}
\phi_n(\theta)= \theta + (\mathbb{Z}\mathbf{D}(\theta_0)\mathbb{Z}')^{-1} \dot{l}_n(\theta).
\end{equation*} 
Then there exists an open ball $B(\theta_0,r)=\lbrace \theta \in \mathbb{R}^k; \| \theta -\theta_0\|< r\rbrace $ (with $ r>0$) such that $\phi_n$ satisfies the Lipschitz condition on $B(\theta_0,r)$ that is,
\begin{equation*}
\| \phi_n(\theta) - \phi_n(\tilde{\theta}\| \leq c\|\theta - \tilde{\theta}\|~\text{for all}~ \theta,\tilde{\theta} \in B(\theta_0,r),~ 0< c < 1.
\end{equation*}
\end{lemma}
Let $\gamma_n$  be defined as $\gamma_n(\theta) = \theta - \phi_n(\theta) = -(\mathbb{Z}\mathbf{D}(\theta)\mathbb{Z}')\dot{l}_n(\theta)$. Then $\gamma_n(\theta_0)$ converges almost surely to 0 as $n \rightarrow \infty$. Note that
\begin{equation*}
\gamma_n(\theta_0) = \left( \frac{1}{n}\ddot{l}_n(\theta_0)\right) ^{-1}\times \left( \frac{1}{n}\dot{l}_n(\theta_0)\right) .
\end{equation*}
The condition C3 implies that, $(\frac{1}{n}\ddot{l}_n(\theta_0))^{-1}$ converges to a matrix $\sum$. Then,
\begin{equation*}
\frac{1}{n}\dot{l}_n(\theta_0) = \frac{1}{n}\mathbb{Z}C(\theta_0) =
\begin{pmatrix}
\frac{1}{n}\sum_{i=1}^n X_{i1}A_i(\theta_0)  \\
\vdots  \\
\frac{1}{n}\sum_{i=1}^n X_{ip}A_i(\theta_0) \\
\frac{1}{n}\sum_{i=1}^n W_{i1}B_i(\theta_0)\\
\vdots \\
\frac{1}{n}\sum_{i=1}^n W_{iq}B_i(\theta_0)
 \end{pmatrix}
\end{equation*}
converges almost surely to $0$ as $n \rightarrow \infty$. Indeed: for all $i =1,\ldots, n$ and $ j= 1,\ldots, p$ we have
\begin{equation*}
\mathbb{E}\left[ X_{ij}A_i(\theta_0)\right]  = \mathbb{E}\left[ \mathbb{E}(X_{ij}A_i(\theta_0)|\mathbf X_i,\mathbf W_i)\right]  = \mathbb{E}\left[ X_{ij}\mathbb{E}(A_i(\theta_0)|\mathbf X_i,\mathbf W_i)\right] 
\end{equation*}
where
\begin{equation*}
\begin{split}
\mathbb{E}(A_i(\theta_0)|\mathbf X_i, \mathbf W_i) &= \mathbb{E}[k_i\times\frac{\mathbf w_i f(\mu'\mathbf w_i)h_i(\beta)^{-n_i}F(\beta'\mathbf x_i)}{k_i(\theta)}-(1-k_i)\frac{\mathbf w_i f(\mu'\mathbf w_i)}{1-F(\mu'\mathbf w_i)}|\mathbf X_i,\mathbf W_i] \\
  &=\frac{w_if(w_i\mu)h_i(\beta)^{-n_i}F(x_i\beta)}{k_i(\theta)}\mathbb{E}[k_i|\mathbf X_i, \mathbf W_i]- \mathbb{E}[(1-k_i)|\mathbf X_i, \mathbf W_i]\frac{w_if(w_i\mu)}{1-F(w_i\mu)}|\mathbf X_i, \mathbf W_i] 
\end{split}
\end{equation*}
We have
\begin{equation*}
\begin{split}
E(J_i|\mathbf X_i, \mathbf W_i)  &= \mathbb{P} (Y_i =0 |\mathbf X_i, \mathbf W_i)\\
                &= p_i + (1- \pi_i)^{m_i}(1-p_i)\\
               &= F(\beta'\mathbf x_i) +(1-F(\mu'\mathbf w_i)^{-n_i}(1-F(\beta'\mathbf x_i)\\
                &= F(\beta'\mathbf x_i) +(1-F(\mu'\mathbf w_i)^{-n_i}h_i(\beta)
\end{split}
\end{equation*}
with
\begin{equation*}
\mathbb{E}[(1-J_i)\mathbb Z_i|\mathbf X_i, \mathbf W_i] = n_i(1-p_i) \pi_i =n_i F(\mu'\mathbf w_i)h_i(\beta).
\end{equation*}
Then, we obtain $\mathbb{E}(A_i(\theta_0)|\mathbf X_i, \mathbf W_i) = 0 $. Similarly, we show that $\mathbb{E}(B_i(\theta_0)|\mathbf X_i, \mathbf W_i)=0$. Moreover, we have
\begin{equation*}
\begin{split}
 \mathbb{V}ar(W_{iq}B_i(\theta_0)) &= \mathbb{E}(\mathbb{V}ar(W_{iq}B_i(\theta_0))|\mathbf X_i, \mathbf W_i)] + \mathbb{V}ar[\mathbb{E}(W_{iq}B_i(\theta_0))|\mathbf X_i, \mathbf W_i]\\
                          &= \mathbb{E}^2_{iq}[\mathbb{V}ar(B_i(\theta_0)|\mathbf X_i, \mathbf W_i)]\\
                          &=\mathbb{E}^2_{iq}(\frac{n_i \mathbf{x}_i f(\beta'\mathbf x_i)}{F(\mu'\mathbf w_i)k_i(\theta)}+\frac{\mathbf x_i f'(\beta'\mathbf x_i)\mathbf z_i-n_i F(\beta'\mathbf x_i))}{F(\beta'\mathbf x_i)h_i(\beta)})]^2 \mathbb{V}ar(k_i|\mathbf X_i, \mathbf W_i)]\\
                          &\leq \mathbb{E}^2_{ik}(\frac{n_i \mathbf x_i f(\beta'\mathbf x_i)}{F(\mu'\mathbf w_i)k_i(\theta)}+\frac{\mathbf x_i f'(\beta'\mathbf x_i)\mathbf z_i-n_i F(\beta'\mathbf x_i))}{F(\beta'\mathbf x_i)h_i(\beta)})^2 = c.
\end{split}
\end{equation*}
The conditions C1, C2 and C3 implie that $ c <\infty$, then 
\begin{equation*}
\sum_{i=1}^n \frac{\mathbb{V}ar(A_i(\theta_0)}{i^2}\leq c \sum_{i=1}^n \frac{1}{i^2}\leq \infty.
\end{equation*}
Kolmogorov's strong law of large numbers ensures that, as $n \rightarrow \infty$
\begin{equation*}
\frac{1}{n}\sum_{i=1}^n\lbrace X_{ij}A_i(\theta_0)- E(X_{ij}A_i(\theta_0))\rbrace = \frac{1}{n}\sum_{i=1}^n A_i(\theta_0)X_{ij} ~ \rightarrow ~0
\end{equation*}
Finaly, as $n \rightarrow \infty$ $\frac{1}{n}\dot{l}_n(\theta_0)$ and $\gamma_n(\theta_0)$ converge almost surely to $0$.
By using similar arguments we proove that, $\frac{1}{n}\sum_{i=1}^n W_{iq}B_i(\theta_0)$ converges almost surely to $0$ as $n \rightarrow \infty$.\\
Let $\epsilon > 0$, the almost surely convergence of $\gamma_n(\theta_0)$ implies then for almost every $\omega \in \Omega$, there exists an integer value $n(\epsilon,\omega)$ such that fo any $ n\geq n(\epsilon,\omega)$, $\|\gamma_n(\theta_0)\|\leq \epsilon $ or equivalently, $0 \in B(\gamma_n(\theta_0),\epsilon)$. In particulary, let $\epsilon = (1-c)s$ with $ 0<c<1$ such as in Lemma \ref{lem1}. Since $\gamma_n$ satisfies the Lipschitz condition, the Lemma \ref{lem1} ensures that there exists an element of $B(\theta_0,s)$ (let denote this element by $\widehat \theta_n$) such that $\gamma_n(\widehat \theta_n) =0$, that is
\begin{equation*}
(\mathbb{Z}\mathbf{D}(\theta_0)\mathbb{Z}')^{-1}\times\dot{l}_n(\widehat \theta_n) = 0.
\end{equation*} 
The condition C3 implies that $\dot{l}_n(\widehat \theta_n) =0$ and that $\widehat \theta_n$ is the unique maximizer
of $l_n$. To summarize, we have shown that for almost every $\omega\in \Omega$ and for every $s > 0$, there exists an integer value $n(s,\omega)$ such that if $n\geq n(s,\omega)$ then the maximum likelihood estimator $\widehat \theta_n$ exists, and $\|\widehat \theta_n - \theta_0 \|\leq s$ (that is, $\widehat \theta_n$ converges almost surely to $\theta_0$).\\
\normalfont
which concludes the proof.\\\\
We now turn to the convergence in distribution of the proposed estimator, which is stated by the following theorem:
\begin{theorem}[Asymptotic normality]\label{th2}
Assume that the conditions \textbf{C1}-\textbf{C5} hold and that $\widehat \theta_n$ converges almost surely to $\theta_0$. Let $\widehat \sum_n = \mathbb{Z}\mathbf{D}(\widehat \theta_n)\mathbb{Z}'$ and $I_k$ denote the identity matrix of order $k$. Then $\widehat{\Sigma}_{n}^{\frac{1}{2}}(\widehat \theta_n - \theta_0)$ converges in distribution to the Gaussian vector $\mathcal{N}(0,I_k)$.
\end{theorem}
\noindent
\textbf{Proof of Theorem \ref{th2}} A Taylor expansion of the score function is as
\begin{equation*}
0 = \dot{l}_n(\widehat \theta_n) = \dot{l}_n(\theta_0) + \ddot{l}_n(\tilde{\theta}_n)(\widehat \theta_n - \theta_0)
\end{equation*}
where $\tilde{\theta}_n$ lies between $\widehat \theta_n$ and $\theta_0$, and thus $\dot{l}_n(\theta_0) = -\ddot{l}_n(\tilde{\theta}_n)(\widehat \theta_n - \theta_0)$. Let $\tilde{\Sigma}_n = -\ddot{l_n}(\tilde{\theta}_n) = \mathbb{Z}\mathbf{D}(\tilde{\theta}_n)\mathbb{Z}' $ and $\Sigma_{n,0} = \mathbb{Z}\mathbf{D}(\theta_0)\mathbb{Z}'$. Now
\begin{equation}\label{eq1th2}
\widehat{\Sigma}_{n}^{\frac{1}{2}}(\widehat \theta_n - \theta_0) = \left[ \widehat{\Sigma}_{n}^{\frac{1}{2}}\tilde{\Sigma}_{n}^{\frac{1}{2}}\right] \left[ \tilde{\Sigma}_n^{\frac{1}{2}}\Sigma_{n,0}^{\frac{1}{2}}\right] \Sigma_{n,0}^{\frac{1}{2}}\left[ \tilde{\Sigma}_n(\widehat \theta_n -\theta_0)\right]. 
\end{equation}
The two terms in brackets in (\ref{eq1th2}) converge almost surely to $I_k$. To see this, we show for example that
$\|\tilde{\Sigma}_{n}^{-\frac{1}{2}}\Sigma_{n,0}^{\frac{1}{2}} - I_k \|\longrightarrow 0$ almost surely as $n\rightarrow \infty$. First, note that 
\begin{equation}\label{eq2th2}
\begin{split}
\|\tilde{\Sigma}_{n}^{-\frac{1}{2}}\Sigma_{n,0}^{\frac{1}{2}} - I_k \| &= \| \tilde{\Sigma}_{n}^{-\frac{1}{2}}\left( \Sigma_{n,0}^{\frac{1}{2}} - \tilde{\Sigma}_{n}^{\frac{1}{2}}\right) \|\\                      
 &\leq \|\tilde{\Sigma}_{n}^{\frac{1}{2}}\| \| \Sigma_{n,0}^{\frac{1}{2}} - \tilde{\Sigma}_{n}^{\frac{1}{2}} \| \\
 &\leq \Lambda_{n}^{\frac{1}{2}}\| \tilde{\Sigma}_{n}^{\frac{1}{2}}  \| \|\Lambda_{n}^{-\frac{1}{2}}\left( \Sigma_{n,0}^{\frac{1}{2}} - \tilde{\Sigma}_{n}^{\frac{1}{2}}\right) \|
\end{split}
\end{equation}
and
\begin{equation*}
\Lambda_{n}^{-1}\| \Sigma_{n,0} -\tilde{\Sigma}_n \| = \Lambda_{n}^{-1}\|\mathbb{Z}(\mathbf{D}(\theta_0) - \mathbf{D}(\tilde{\theta}_n))\mathbb{Z}'\|.
\end{equation*}
Note also that $\widehat \theta_n$ converges almost surely to $\theta_0$ (that is, for every $w\in \tilde{\Omega}$, where $\tilde{\Omega}\subset \Omega$ and $\mathbb{P}(\tilde{\Omega})=1$). Let $w\in \tilde{\Omega}$. By the same arguments as in the proof of Lemma \ref{lem1}, for every $\epsilon > 0$, there exists a positive $n(\epsilon,w)\in \mathbb{N}$ such that if $n\geq n(\epsilon,w)$, then $\Lambda_{n}^{-1}\|\mathbb{Z}(\mathbf{D}(\theta_0) - \mathbf{D}(\tilde{\theta}_n))\mathbb{Z}'\|< \epsilon $. Hence $\Lambda_{n}^{-1}\|\mathbb{Z}(\mathbf{D}(\theta_0) - \mathbf{D}(\tilde{\theta}_n))\mathbb{Z}'\|$ converges almost surely to $0$. By continuity of the map $x\rightarrow x^{\frac{1}{2}}$, $\|\Lambda_{n}^{-\frac{1}{2}}\left( \Sigma_{n,0}^{\frac{1}{2}} - \tilde{\Sigma}_{n}^{\frac{1}{2}}\right) \|$ converges also almost surely to $0$. Moreover, for $n$ sufficiently large, there exists a positive constant $A < \infty$, such that almost surely, $\|\tilde{\Sigma}_{n}^{-\frac{1}{2}}\|\leq A\frac{\Lambda_{n}^{\frac{1}{2}}}{\lambda_{n}^{1/2}} < A c_{2}^{\frac{1}{2}}$. It follows from (\ref{eq2th2}) and the condition C3 that $\|(\Sigma_{n,0}^{\frac{1}{2}}\tilde{\Sigma}_{n}^{-\frac{1}{2}}) -I_k \|$ converges almost surely to $0$. 
It remains for us to show that $\Sigma_{n,0}^{-\frac{1}{2}}(\tilde{\Sigma}_n(\widehat \theta_n - \theta_0))$ converges in distribution to $\mathcal{N}(0,I_k)$. Note that $\Sigma_{n,0}^{-\frac{1}{2}}(\tilde{\Sigma}_{n}(\widehat \theta_n - \theta_0)) = \Sigma_{n,0}^{-\frac{1}{2}}\sum_{j=1}^{2n}\mathbb{Z}_{\bullet j}C_j(\theta_0)$. This convergence holds if we can check the following conditions:
\begin{itemize}
\item[i)] $\max_{1\leq j\leq 2n}\mathbb{Z}'_{\bullet j}(\mathbb{Z}\mathbb{Z}')^{-1}\mathbb{Z}_{\bullet j}\rightarrow 0$ if $n\rightarrow \infty$. Condition i) follows by noting that 
\begin{equation*}
0\leq \max_{1\leq j\leq 2n}\mathbb{Z}'_{\bullet j}(\mathbb{Z}\mathbb{Z}')^{-1}\mathbb{Z}_{\bullet j} \leq \max_{1\leq j\leq 2n}\|\mathbb{Z}_{\bullet j}\|^2\|(\mathbb{Z}\mathbb{Z}')^{-1}\| = \max_{1\leq j\leq 2n}\frac{\| \mathbb{Z}_{\bullet j} \|^2}{\tilde{\lambda}_n}.
\end{equation*}
and that $\| \mathbb{Z}_{\bullet j}\|^2$ is bounded above, by C1 and C2. Moreover, $\frac{1}{\tilde{\lambda}_n}\rightarrow 0 $ as $ n\rightarrow \infty$.
\item[ii)] $\sup_{1\leq j \leq 2n}\mathbb{E}\left[ C_{j}^{2}(\theta_0) \mathbf{1}_{\lbrace |C_j(\theta_0)|>c\rbrace}\right]  \rightarrow 0 $ as $n \rightarrow \infty$. Condition ii) follows by noting that the components $C_j(\theta_0)$ of $C$ are bounded under C1, C2 and C3. Finally, for every $i = 1, \ldots, 2n$, $ \mathbb{E}(C_{j}^{2}(\theta_0)) = Var(C_j(\theta_0)) $ since $ \mathbb{E}(C_j(\theta_0)) = 0$. If $ u\in \lbrace n+1,\ldots,2n\rbrace$, $C_j(\theta_0) = B_{h}(\theta_0)$ with $ u = j-n$ then
\begin{equation*}
\begin{split}
\mathbb{V}ar(C_j(\theta_0)) &= \mathbb{V}ar(B_{h}(\theta_0)) \\
&= \mathbb{E}(\mathbb{V}ar(B_{h}(\theta_0))|\mathbf X_{u}, \mathbf W_{u})) + \mathbb{V}ar(\mathbb{E}(B_{h}(\theta_0)|\mathbf X_{u}, \mathbf W_{u})) \\
&= \mathbb{E}(\mathbb{V}ar(B_{h}(\theta_0))| \mathbf X_{u}, \mathbf W_{u}).
\end{split}
\end{equation*}
Finally,
\begin{equation*}
\begin{split}
Var(B_{h}(\theta_0)| \mathbf X_{h}, \mathbf W_{h}) &= \left(  \frac{-k_i n_i \mathbf x_i f(\beta'\mathbf x_i)}{F(\mu'\mathbf w_i)k_i(\theta)} + (1-k_i)\frac{\mathbf x_i f(\beta'\mathbf x_i)[\mathbf z_i- n_i F(\beta'\mathbf x_i)]}{F(\beta'\mathbf x_i)h_{i}}(\beta)\right)^{2} \\
&\times Var(J_{u}| \mathbf X_{u}, \mathbf W_{u})\\
                 &= \left( \frac{-k_i n_i \mathbf x_i f(\beta'\mathbf x_i)}{F(\mu'\mathbf w_i)k_i(\theta)} + (1-k_i)\frac{\mathbf x_i f(\beta'\mathbf x_i)[\mathbf z_i-n_i F(\beta'\mathbf x_i)]}{F(\beta'\mathbf x_i)h_{i}(\beta)}\right)^{2}\\
                 &\times \mathbb{P}(Y_{u}= 0|\mathbf X_{u},\mathbf W_{u})(1-\mathbb{P}(Y_{u}=0|\mathbf X_{u},\mathbf W_{u}))\\
                 &= \left( \frac{-k_i n_i \mathbf x_i f(\beta'\mathbf x_i)}{F(\mu'\mathbf w_i)k_i(\theta)} + (1-k_i)\frac{\mathbf x_i f(\beta'\mathbf x_i)[\mathbf z_i-n_i F(\beta'\mathbf x_i)]}{F(\beta'\mathbf x_i)h_{i}(\beta)}\right)^{2} \\
                 &\times \left( p_{u}+(1-\pi_{u})^{m_u}(1-p_{j}))(1-p_{u})\times (1-(1-\pi_{u})^{m_u}\right).              
\end{split}
\end{equation*}

Therefore $\mathbb{V}ar(B_{u}(\theta_0)| \mathbf X_{u}, \mathbf W_{u}) > 0$ and for all $j \in \lbrace 1,\ldots,n\rbrace $ under conditions C1, C2 and C3, $\mathbb{V}ar(C_j(\theta_0))> 0$.
\item[iii)] $\inf_{1\leq j\leq 2n}\mathbb{E}(C_{j}^{2}(\theta_0))>0$. Condition iii) follows by using similar arguments,  that is $\mathbb{V}ar(C_j(\theta_0))> 0$ for all $j = 1,\ldots, 2n$ 
\end{itemize}
To summarize, we have proved that   $\Sigma_{n,0}^{-\frac{1}{2}}(\tilde{\Sigma}_n(\widehat \theta_n - \theta_0))$ converges in distribution to $\mathcal{N}(0,I_k)$. This result, combined with Slutsky’s theorem and equation (\ref{eq1th2}), implies that $\widehat \Sigma_{n}^{\frac{1}{2}}(\widehat \theta_n - \theta_0)$ converges in distribution to $\mathcal{N}(0,I_k)$.\\
\normalfont
which concludes the proof.
\section{Discussion and perspectives}
\label{sec:discussion}
In this paper, we have proposed a procedure of maximum likelihood estimation in the Zero-inflated Binomial regression model by using probit link function for both coponents of the model (\ref{sec1eq1}). We establish the asymptotic properties (existence, consistency and asymptotic normality of the proposed maximum likelihood estimator).\\
Several open problems still deserve attention. First, some simulations to study the proposed of the estimator in finite sample and application on real data. It is sometimes of interest to make inference about the probability of event $\mathbb P(Y=y_i| \mathbf W)$ across the whole range of the predictors $\mathbf W$. The calculation of simultaneous confidence bands for the probabilities $\{p(\mathbf x), \mathbf x\in \mathcal X\}$ thus constitutes another issue of (both methodological and practical) interest. Another issue of interest deals with the inference in the ZIB regression model with probit link function, in a high-dimensional setting. We have established the theoretical properties of our estimator in a low-dimensional setting that is, when a small number of potential predictors are involved (this problem arises, for example, in genetic studies where high-dimensional data are generated using microarray technologies).

\bibliographystyle{Chicago}

\end{document}